\documentclass[11pt]{amsart}
\usepackage{amsmath,amssymb,amsthm,latexsym,mathrsfs}
\usepackage{indentfirst}
 \DeclareMathOperator{\DIV}{div}
  \DeclareMathOperator{\grad}{grad}
 \DeclareMathOperator{\trace}{tr}

 \newcommand{\dif}{\mathrm{d}}
 \newcommand{\ROM}[1]{\mathrm{\uppercase\expandafter{\romannumeral#1}}}
  \theoremstyle{definition}
 \newtheorem{theorem}{Theorem}[section]
 \newtheorem{lemma}[theorem]{Lemma}

 \newtheorem{remark}[theorem]{Remark}
 \newtheorem{proposition}[theorem]{Proposition}
 
 \numberwithin{equation}{section}
   \makeatletter
 \renewcommand{\theequation}{\@arabic\c@section.\arabic{equation}}
\renewcommand \thesection {\S\@arabic\c@section}
\renewcommand \thetheorem {\@arabic\c@section.\@arabic\c@theorem}
\makeatother
\title[Stability of hypersurfaces]{\textbf{Stability of hypersurfaces with constant $r$-th anisotropic mean curvature}}
\author[Y. J. He]{Yijun He}\address{School of Mathematical Sciences,
Shanxi University, Taiyuan 030006, P. R. China.}
\thanks {The first author was partially supported by Youth Science Foundation
of Shanxi Province, China (Grant No. 2006021001).}
\email{heyijun@sxu.edu.cn}
 \author[H. Li]{Haizhong
Li}\address{Department of Mathematical Sciences, Tsinghua
University, Beijing 100084, P. R.
China.}\email{hli@math.tsinghua.edu.cn}
\thanks {The second author was partially supported by the grant No. 10531090 of the NSFC and by
SRFDP.}
\date{}
\begin{document}
\maketitle
\begin{abstract}
 Given a positive function $F$ on $S^n$ which satisfies a convexity
condition, we define the $r$-th anisotropic mean curvature function
$H^F_r$ for hypersurfaces in $\mathbb{R}^{n+1}$ which is a
generalization of the usual $r$-th mean curvature function. Let
$X:M\to \mathbb{R}^{n+1}$ be an $n$-dimensional closed hypersurface
with $H^F_{r+1}=$constant, for some $r$ with $0\leq r\leq n-1$,
which is a critical point for a variational problem. We show that
$X(M)$ is stable if and only if $X(M)$ is the Wulff shape.
\end{abstract}
\section{Introduction}
Let $F\colon S^n\to\mathbb{R}^+$ be a smooth function which
satisfies the following convexity condition:
\begin{equation}\label{1}
(D^2F+F1)_x>0,\quad\forall x\in S^n,
\end{equation}
where $S^n$ denotes the standard unit sphere in $\mathbb{R}^{n+1}$,
$D^2F$ denotes the intrinsic Hessian of $F$ on $S^n$ and 1 denotes
the identity on $T_x S^n$, $>0$ means that the matrix is positive
definite. We consider the map
\begin{equation}\label{2}
\begin{array}{c}
\phi\colon S^n\to\mathbb{R}^{n+1},\\
 x\mapsto F(x)x+(\grad_{S^n}F)_x,
 \end{array}
\end{equation}
its image $W_F=\phi(S^n)$ is a smooth, convex hypersurface in
$\mathbb{R}^{n+1}$ called the Wulff shape of $F$ (see
\cite{clarenz}, \cite{HeL1}, \cite{HeL2}, \cite{K-Pa},
\cite{palmer}, \cite{Taylor}, \cite{winklmann}). We note when
$F\equiv1$, $W_F$ is just the sphere $S^n$.

Now let $X\colon M\to\mathbb{R}^{n+1}$ be a smooth immersion of a
closed, orientable hypersurface. Let $\nu\colon M\to S^n$ denotes
its Gauss map, that is, $\nu$ is the unit inner normal vector of
$M$.

Let $A_F=D^2F+F1$, $S_F=-\dif(\phi\circ\nu)=-A_F\circ\dif\nu$. $S_F$
is called the $F$-Weingarten operator, and the eigenvalues of $S_F$
are called anisotropic principal curvatures. Let $\sigma_r$ be the
elementary symmetric functions of the anisotropic principal
curvatures $\lambda_1, \lambda_2, \cdots, \lambda_n$:
\begin{equation}\label{3}
\sigma_r=\sum_{i_1<\cdots<i_r}\lambda_{i_1}\cdots\lambda_{i_r}\quad
(1\leq r\leq n).
\end{equation}
 We set $\sigma_0=1$. The $r$-th
anisotropic mean curvature $H^F_r$ is defined by
$H^F_r=\sigma_r/C^r_n$, also see Reilly \cite{reilly2}.

For each $r$, $0\leq r\leq n-1$, we set
\begin{equation}\label{4}\mathscr{A}_{r, F}=\int_MF(\nu)\sigma_r\dif A_X.
\end{equation}

The algebraic $(n+1)$-volume enclosed by $M$ is given by
\begin{equation}\label{5}
  V=\frac1{n+1}\int_M\langle X, \nu\rangle\dif A_X.
\end{equation}

We consider those hypersurfaces which are critical points of
$\mathscr{A}_{r, F}$ restricted to those hypersurfaces enclosing a
fixed volume $V$. By a standard argument involving Lagrange
multipliers, this means we are considering critical points of the
functional
\begin{equation}\label{6}
\mathscr{F}_{r, F; \Lambda}=\mathscr{A}_{r, F}+\Lambda V(X),
\end{equation}
where $\Lambda$ is a constant. We will show the Euler-Lagrange
equation of $\mathscr{F}_{r, F; \Lambda}$ is:
\begin{equation}\label{7}
  (r+1)\sigma_{r+1}-\Lambda=0.
\end{equation}
 So the critical points are just hypersurfaces with $H^F_{r+1}=\mbox{const}$.

 If $F\equiv1$, then the function $\mathscr{A}_{r, F}$ is just the
 functional $\mathscr{A}_r=\int_MS_r\dif A_X$ which was studied by Alencar, do Carmo and Rosenberg in  \cite{AdR}, where $H_r=S_r/C_n^r$
 is the usual $r$-th mean curvature. For such a variational problem, they call a critical immersion $X$
 of the  functional $\mathscr{A}_r$ (that is, a hypersurface with $H_{r+1}=\mbox{constant}$)
 stable if and only if the second
variation of $\mathscr{A}_r$ is non-negative for all variations of
$X$ preserving the enclosed $(n+1)$-volume $V$. They proved:
\begin{theorem}\label{Adr} (\cite{AdR})
  Suppose $0\leq r\leq n-1$. Let $X\colon M\to\mathbb{R}^{n+1}$ be a closed hypersurface with $H_{r+1}=\mbox{constant}$. Then $X$
  is stable if and only if $X(M)$ is a round sphere.
\end{theorem}

Analogously, we call a critical immersion $X$ of the functional
$\mathscr{A}_{r, F}$ stable if and only if the second variation of
$\mathscr{A}_{r, F}$ (or equivalently of $\mathscr{F}_{r, F;
\Lambda}$) is non-negative for all variations of $X$ preserving the
enclosed $(n+1)$-volume $V$.

In \cite{palmer}, Palmer proved the following theorem (also see
Winklmann \cite{winklmann}):
\begin{theorem} (\cite{palmer})
 Let $X\colon M\to\mathbb{R}^{n+1}$ be a closed hypersurface with $H^F_1=$constant. Then $X$
  is stable if and only if, up to translations and homotheties, $X(M)$ is the Wulff
 shape.
\end{theorem}

In this paper, we prove the following theorem:
\begin{theorem}\label{thm1.2}
  Suppose $0\leq r\leq n-1$. Let $X\colon M\to\mathbb{R}^{n+1}$ be a closed hypersurface with $H^F_{r+1}=$constant.
  Then, $X$ is stable if and only if,
   up to translations and homotheties, $X(M)$ is the Wulff
 shape.
\end{theorem}
\begin{remark}
  In  the case $F\equiv1$, Theorem \ref{thm1.2}  becomes Theorem \ref{Adr}.
  Theorem \ref{thm1.2} gives an affirmative
answer to the problem proposed in \cite{HeL2}.
\end{remark}
\section{Preliminaries}

Let $X\colon M\to R^{n+1}$ be a smooth closed, oriented hypersurface
with Gauss map $\nu\colon M\to S^n$, that is, $\nu$ is the unit
inner normal vector field. Let $X_t$ be a variation of $X$, and
$\nu_t\colon M\to S^n$ be the Gauss map of $X_t$. We define
\begin{equation}\label{8}
\psi=\langle\frac{\dif X_t}{\dif t},\nu_t\rangle,\quad
\xi=(\frac{\dif X_t}{\dif t})^{\top},
\end{equation}
where $\top$ represents the tangent component and $\psi$, $\xi$ are
dependent of $t$. The corresponding first variation of the unit
normal vector is given by (see \cite{K-Pa}, \cite{palmer},
\cite{winklmann})
\begin{equation}\label{9}
\nu_t'=-\grad\psi+\dif\nu_t(\xi),
\end{equation}
the first variation of the volume element is (see \cite{BC},
\cite{CL} or \cite{HL})
\begin{equation}\label{10}
\partial_t\dif A_{X_t}=(\DIV\xi-nH\psi)\dif A_{X_t},
\end{equation}
and the first variation of the volume $V$ is
 \begin{equation}\label{11}
 V'(t)=\int_M\psi dA_{X_t},
\end{equation}
where
$\grad$, $\DIV$, $H$ represents the gradients, the divergence, the
mean curvature with respect to $X_t$ respectively.

Let $\{E_1,\cdots,E_n\}$ be a local orthogonal frame on $S^n$, let
$e_i=e_i(t)=E_i\circ\nu_t$, where $i=1,\cdots,n$ and $\nu_t$ is the
Gauss map of $X_t$, then $\{e_1,\cdots,e_n\}$ is a local orthogonal
frame of $X_t\colon M\to\mathbb{R}^{n+1}$.

The structure equation of $x\colon S^n\to\mathbb{R}^{n+1}$ is:
\begin{equation}\label{12}\left\{
\begin{array}
  {l}
  \dif x=\sum_i\theta_iE_i\\
  \dif E_i=\sum_j\theta_{ij}E_j-\theta_ix\\
  \dif\theta_i=\sum_j\theta_{ij}\wedge\theta_j\\
  \dif\theta_{ij}-\sum_k\theta_{ik}\wedge\theta_{kj}=-\frac12\sum_{kl}{\tilde R}_{ijkl}\theta_k\wedge\theta_l=-\theta_i\wedge\theta_j
\end{array}\right.
\end{equation}
where $\theta_{ij}+\theta_{ji}=0$ and
 ${\tilde R}_{ijkl}=\delta_{ik}\delta_{jl}-\delta_{il}\delta_{jk}$.

The structure equation of $X_t$ is (see \cite{L1}, \cite{L2}):
\begin{equation}\label{14}\left\{
\begin{array}
  {l}
  \dif X_t=\sum_i\omega_ie_i\\
  \dif \nu_t=-\sum_{ij}h_{ij}\omega_je_i\\
  \dif
  e_i=\sum_j\omega_{ij}e_j+\sum_jh_{ij}\omega_j\nu_t\\
  \dif\omega_i=\sum_j\omega_{ij}\wedge\omega_j\\
  \dif\omega_{ij}-\sum_k\omega_{ik}\wedge\omega_{kj}=
  -\frac12\sum_{kl}R_{ijkl}\theta_k\wedge\theta_l
\end{array}\right.
\end{equation}
where $\omega_{ij}+\omega_{ji}=0$, $R_{ijkl}+R_{ijlk}=0$, and $R_{ijkl}$ are the components of the Riemannian curvature tensor of $X_t(M)$ with
respect to the induced metric $dX_t\cdot dX_t$. Here we have omitted the variable $t$ for some geometric quantities.

From $\dif e_i=\dif(E_i\circ\nu_t)=\nu_t^*\dif E_i=\sum_j\nu_t^*\theta_{ij}e_j-\nu_t^*\theta_i\nu_t$, we get
\begin{equation}\label{15}\left\{
\begin{array}
  {l}
  \omega_{ij}=\nu_t^*\theta_{ij}\\
  \nu_t^*\theta_i=-\sum_jh_{ij}\omega_j,
\end{array}\right.
\end{equation}
where $\omega_{ij}+\omega_{ji}=0$, $h_{ij}=h_{ji}$.

 Let $F\colon S^n\to\mathbb{R}^+$ be a smooth function, we denote the coefficients of covariant differential of $F$, $\grad_{S^n}F$
 with respect to $\{E_i\}_{i=1,\cdots,n}$ by $F_i, F_{ij}$
respectively.

From (\ref{15}), $\dif(F(\nu_t))=\nu_t^*\dif
F=\nu_t^*(\sum_iF_i\theta_i)=-\sum_{ij}(F_i\circ\nu_t)h_{ij}\omega_j$,
thus
\begin{equation}
  \label{gradf}
  \grad
  (F(\nu_t))=-\sum_{ij}(F_i\circ\nu_t)h_{ij}e_j=\dif\nu_t(\grad_{S^n}F).
\end{equation}

Through a direct calculation, we easily get
\begin{equation}
  \label{dphi}
  \dif\phi=(D^2F+F1)\circ \dif x=\sum_{ij}A_{ij}\theta_iE_j,
\end{equation}
where $A_{ij}$ is the coefficient of $A_F$, that is,
$A_{ij}=F_{ij}+F\delta_{ij}$.

Taking exterior differential of (\ref{dphi}) and using (\ref{12}) we
get
\begin{equation}\label{17}
A_{ijk}=A_{jik}=A_{ikj},
\end{equation}
where $A_{ijk}$ denotes coefficient of the covariant differential of $A_F$ on $S^n$.

We define $(A_{ij}\circ\nu_t)_k$ by
\begin{equation}\label{20}
\dif(A_{ij}\circ\nu_t)+\sum(A_{kj}\circ\nu_t)\omega_{ki}+\sum_k(A_{ik}\circ\nu_t)\omega_{kj}=
\sum_k(A_{ij}\circ\nu_t)_k\omega_k.
\end{equation}

By a direct calculation using (\ref{15}) and (\ref{20}), we have
\begin{equation}\label{21}
(A_{ij}\circ\nu_t)_k=-\sum_lh_{kl}A_{ijl}\circ\nu_t.
\end{equation}

We define $L_{ij}$ by
 \begin{equation}\label{22}(\frac{\dif e_i}{\dif t})^{\top}=-\sum_jL_{ij}e_j,
\end{equation}
 where $\top$ denote the tangent component,
 then $L_{ij}=-L_{ji}$ and we have (see \cite{BC}, \cite{CL} or \cite{HL})
 \begin{equation}
 \label{23}
   h'_{ij}=\psi_{ij}+\sum_k\{h_{ijk}\xi_k+\psi
   h_{ik}h_{jk}+h_{ik}L_{kj}+h_{jk}L_{ki}\}.
\end{equation}

 Let $s_{ij}=\sum_kA_{ik}h_{kj}$, then from (\ref{15}) and (\ref{dphi}), we have
\begin{equation}
  \label{dphinu}
  \dif(\phi\circ\nu_t)=\nu_t^*\dif\phi=-\sum_{ij}s_{ij}\omega_je_i.
\end{equation}
We define $S_F$ by $S_F=-\dif(\phi\circ\nu)=-A_F\circ\dif\nu$, then
we have $S_F(e_j)=\sum_is_{ij}e_i$. We call $S_F$ to be F-Weingarten
operator. From the positive definite of $(A_{ij})$ and the symmetry
of $(h_{ij})$, we know the eigenvalues of $(s_{ij})$ are all real.
we call them anisotropic principal curvatures, and denote them by
$\lambda_1, \cdots, \lambda_n$.

Taking exterior differential of (\ref{dphinu}) and using (\ref{14})
we get
\begin{equation}\label{sijk}
s_{ijk}=s_{ikj},
\end{equation}
where $s_{ijk}$ denotes coefficient of the covariant differential of
$S_F$.

We have $n$ invariants, the elementary symmetric function $\sigma_r$ of the anisotropic principal curvatures:
\begin{equation}\label{24}
  \sigma_r=\sum_{i_1<\cdots i_r}\lambda_{i_1}\cdots\lambda_{i_n}
  \quad (1\leq r\leq n).
\end{equation}
For convenience, we set $\sigma_0=1$ and $\sigma_{n+1}=0$. The
$r$-th anisotropic mean curvature $H^F_r$ is defined by
\begin{equation}\label{25}
H^F_r=\sigma_r/C_n^r,\quad C^r_n=\frac{n!}{r!(n-r)!}.
\end{equation}

We have by use of (\ref{9}) and (\ref{14})
 \begin{equation}\label{26}
 \begin{array}{rl}
 &\sum_{ij}\dfrac{\dif((A_{ij}E_i\otimes E_j)\circ\nu_t)}{\dif t}=\sum_{ij}\langle (D(A_{ij}E_i\otimes E_j))_{{\nu}_t}, \nu_t'\rangle\\
 =& -\sum_{ijk}A_{ijk}(\psi_k+\sum_lh_{kl}\xi_l)e_i\otimes e_j,\end{array}
\end{equation}
where $D$ is the Levi-Civita connection on $S^n$.

 On the other hand, we have
\begin{equation}\label{27}
 \sum_{ij}\frac{\dif((A_{ij}E_i\otimes E_j)\circ\nu_t)}{\dif t}=
 \sum_{ij}\{A_{ij}'e_i\otimes e_j+A_{ij}(\frac{\dif e_i}{\dif t})^{\top}\otimes e_j+A_{ij}e_i\otimes (\frac{\dif e_j}{\dif t})^{\top}\}.
\end{equation}

 By use of (\ref{22}), we get from (\ref{26}) and (\ref{27})
\begin{equation}\label{28}
   \frac{d(A_{ij}\circ \nu_t)}{dt}=A'_{ij}(t)=\sum_k\{-A_{ijk}\psi_k-\sum_lA_{ijk}h_{kl}\xi_l+A_{ik}L_{kj}+A_{jk}L_{ki}\}.
\end{equation}

By (\ref{21}), (\ref{23}), (\ref{28}) and the fact $L_{ij}=-L_{ji}$,
through a direct calculation, we get the following lemma:

\begin{lemma}\label{lem2.1}
$\dfrac{ds_{ij}}{dt}=
s'_{ij}(t)=\sum_k\{(A_{ik}\psi_k)_j+s_{ijk}\xi_k+\psi
   s_{ik}h_{kj}+s_{kj}L_{ki}+s_{ik}L_{kj}\}$.
\end{lemma}

As $M$ is a closed oriented hypersurface,  one can find a point
where all the principal curvatures with respect to $\nu$ are
positive. By the positiveness of $A_F$, all the anisotropic
principal curvatures are positive at this point. Using the results
of G{\aa}rding (\cite{G}), we have the following lemma:
\begin{lemma}\label{lemma3.0}
Let $X\colon M\to\mathbb{R}^{n+1}$ be a closed, oriented
hypersurface. Assume $H^F_{r+1}>0$ holds on every point of $M$, then
 $H^F_k>0$ holds on every point of $M$ for every $k=1, \cdots, r$.
\end{lemma}

 Using the
characteristic polynomial of $S_F$, $\sigma_r$ is defined by
\begin{equation}\label{30}
 \det(tI-S_F)=\sum_{r=0}^n(-1)^r\sigma_rt^{n-r}.
\end{equation}
 So, we have
 \begin{equation}\label{sigmar}
   \sigma_r=\frac{1}{r!}\sum_{i_1,\cdots,i_r;j_1,\cdots,j_r}\delta_{i_1\cdots
   i_r}^{j_1\cdots j_r}s_{i_1j_1}\cdots s_{i_rj_r},
 \end{equation}
  where $\delta_{i_1\cdots
   i_r}^{j_1\cdots j_r}$ is the usual generalized Kronecker symbol,
   i.e., $\delta_{i_1\cdots
   i_r}^{j_1\cdots j_r}$ equals +1 (resp. -1) if $i_1\cdots
   i_r$ are distinct and  $(j_1\cdots j_r)$
   is an even (resp. odd) permutation of $(i_1\cdots i_r)$ and in
   other cases it equals zero.

  We introduce two important operators $P_r$ and $T_r$ by
\begin{equation}\label{pr}
P_r=\sigma_rI-\sigma_{r-1}S_F+\cdots+(-1)^rS_F^r, \quad r=0, 1,
\cdots, n, \end{equation}
\begin{equation}\label{tr}
T_r=P_rA_F, \quad r=0, 1, \cdots, n-1. \end{equation}
 Obviously $P_n=0$ and we have
 \begin{equation}
   \label{pr+1}
   P_r=\sigma_rI-P_{r-1}S_F=\sigma_rI+T_{r-1}\dif\nu, \quad r=1, \cdots, n.
 \end{equation}

  From the symmetry of $A_F$ and
     $\dif\nu$, $S_FA_F$ and $\dif\nu\circ S_F$ are symmetric, so
$T_r=P_rA_F$ and $\dif\nu\circ P_r$ are also symmetric for each $r$.

   \begin{lemma}\label{lemma3.3}
   The matrix of $P_r$ is given by:
   \begin{equation}
 \label{Prij}
 (P_r)_{ij}=\displaystyle{\frac{1}{r!}\sum_{i_1,\cdots,i_r;j_1,\cdots,j_r}
     \delta_{i_1\cdots i_r j}^{j_1\cdots j_r i}}s_{i_1j_1}\cdots s_{i_rj_r}
 \end{equation}
 \end{lemma}
 \begin{proof}
   We prove Lemma \ref{lemma3.3} inductively. For $r=0$, it is easy to check that (\ref{Prij}) is
   true.

Assume (\ref{Prij}) is true for $r=k$, then from (\ref{pr+1}),
$$\begin{array}
  {rcl}
  (P_{k+1})_{ij} &=& \sigma_{k+1}\delta^i_j-\sum_l(P_k)_{il}s_{lj}\\
  &=& \dfrac{1}{(k+1)!}\sum(\delta_{i_1\cdots
       i_{k+1}}^{j_1\cdots
       j_{k+1}}\delta^i_j-\sum_l\delta_{i_1\cdots
       i_{l-1}i_li_{l+1}\cdots i_{k+1}}^{j_1\cdots
       j_{l-1}ij_{l+1}\cdots
       j_{k+1}}\delta_j^{j_l})s_{i_1j_1}\cdots s_{i_{k+1}j_{k+1}}\\
       &=& \dfrac{1}{(k+1)!}\sum\delta_{i_1\cdots
       i_{k+1}j}^{j_1\cdots
       j_{k+1}i}s_{i_1j_1}\cdots s_{i_{k+1}j_{k+1}}.
\end{array}$$

 \end{proof}

   \begin{lemma}\label{lemma3.4}
   For each $r$, we have

   (\romannumeral1). $\sum_j(P_r)_{jij}=0$,

   (\romannumeral2). $\trace(P_rS_F)=(r+1)\sigma_{r+1}$,

   (\romannumeral3). $\trace(P_r)=(n-r)\sigma_r$,

   (\romannumeral4). $\trace(P_rS_F^2)=\sigma_1\sigma_{r+1}-(r+2)\sigma_{r+2}$.

\end{lemma}
\begin{proof}
  We only prove (\romannumeral2), the others are easily obtained from
  (\ref{sigmar}), (\ref{pr+1}) and (\ref{Prij}).

  Noting $(j,j_r)$ is symmetric in $s_{i_1j_1}\cdots
s_{i_rj_rj}$ by (\ref{sijk}) and $(j,j_r)$ is skew symmetric in
$\delta_{i_1\cdots i_r i}^{j_1\cdots j_r j}$  , we have

 $$\sum_{j}(P_r)_{jij}=\displaystyle{\frac{1}{(r-1)!}\sum_{i_1,\cdots,i_r;j_1,\cdots,j_r;j}
     \delta_{i_1\cdots i_r i}^{j_1\cdots j_r j}}s_{i_1j_1}\cdots
     s_{i_rj_rj}=0.
$$
\end{proof}

\begin{remark}
When $F=1$, Lemma \ref{lemma3.4} was a well-known result (for
example, see Barbosa-Colares \cite{BC}, Reilly \cite{reilly}, or
Rosenberg \cite{rosenberg}).
\end{remark}
Since $P_{r-1}S_F$ is symmetric and $L_{ij}$ is anti-symmetric, we
have
\begin{equation}\label{lemma3.5}
  \sum\limits_{i,j,k}(P_{r-1})_{ji}(s_{kj}L_{ki}+s_{ik}L_{kj})=0.
\end{equation}

From (\ref{sijk}), (\ref{pr+1}) and (\romannumeral1) of Lemma
\ref{lemma3.4}, we get
\begin{equation}
  \label{sigmark}
  (\sigma_r)_k=\sum_j(\sigma_r\delta_{jk})_j=\sum_j(P_r)_{jkj}+\sum_{jl}[(P_{r-1})_{jl}s_{lk}]_j=\sum_{ij}(P_{r-1})_{ji}s_{ijk}.
\end{equation}

\section{First and second variation formulas of
$\mathscr{F}_{r, F; \Lambda}$}

Define the operator $L_r: C^\infty(M)\to C^\infty(M)$ as following:
\begin{equation}\label{42}
L_rf=\sum\limits_{i,j}[(T_r)_{ij}f_j]_i. \end{equation}
\begin{lemma}\label{lemma4.1}
$ \dfrac {d\sigma_r}{dt}=\sigma_r'(t)=L_{r-1}\psi+\psi\langle
T_{r-1}\circ\dif\nu_t,
  \dif\nu_t\rangle+\langle\grad\sigma_r, \xi\rangle.$
\end{lemma}
\begin{proof}
  Using (\ref{sigmar}), (\ref{lemma3.5}), (\ref{sigmark}), Lemma \ref{lem2.1}, Lemma \ref{lemma3.3} and (i) of Lemma \ref{lemma3.4},
  we have
   $$\begin{array}{rcl}
    \sigma_r' &=& \dfrac{1}{(r-1)!}\sum_{i_1,\cdots,i_r;j_1,\cdots,j_r}\delta_{i_1\cdots
    i_r}^{j_1\cdots j_r}s_{i_1j_1}\cdots
    s_{i_{r-1}j_{r-1}}s'_{i_rj_r}\\
    &=& \sum_{ijk}(P_{r-1})_{ji}s'_{ij}\\
    &=&\sum_{ijk}(P_{r-1})_{ji}[(A_{ik}\psi_k)_j+\psi
    s_{ik}h_{kj}+s_{ijk}\xi_k+s_{kj}L_{ki}+s_{ik}L_{kj}]\\
    &=& \sum_{ijk}[(P_{r-1})_{ji}A_{ik}\psi_k]_j+\psi\sum_{ijkl}
    (P_{r-1})_{ji}A_{il}h_{lk}h_{kj}+\sum_k(\sigma_r)_k\xi_k\\
    &=&\sum_{jk}[(T_{r-1})_{jk}\psi_k]_j+\psi\sum_{ijk}
    (T_{r-1})_{ji}h_{ik}h_{kj}+\sum_k(\sigma_r)_k\xi_k\\
    &=& L_{r-1}\psi+\psi\langle T_{r-1}\circ\dif\nu_t,
  \dif\nu_t\rangle+\langle\grad\sigma_r, \xi\rangle.
  \end{array}$$
\end{proof}
\begin{lemma}
  \label{newlemma}For each $0\leq r\leq n$, we have
  \begin{equation}
  \DIV(P_r(\grad_{S^n} F)\circ\nu_t)+F(\nu_t)\trace(P_r\circ
  \dif\nu_t)=-(r+1)\sigma_{r+1},\end{equation}
and
  \begin{equation}
  \DIV(P_rX^{\top})+\langle X, \nu_t\rangle\trace(P_r\circ
  \dif\nu_t)=(n-r)\sigma_r.
  \end{equation}
\end{lemma}
\begin{proof}
  From (\ref{14}), (\ref{dphinu}) and Lemma
\ref{lemma3.4},
 $$\begin{array}
  {rcl}
  \DIV(P_r(\grad_{S^n}
F)\circ\nu_t) &=&
\DIV(P_r(\phi\circ\nu_t)^{\top})\\
&=& \sum_{ij}((P_r)_{ji}\langle\phi\circ\nu_t,
e_i\rangle)_j\\
&=& -\sum_{ij}(P_r)_{ji}s_{ij}+F(\nu_t)\sum_{ij}(P_r)_{ji}h_{ij}\\
&=& -\trace(P_rS_F)-F(\nu_t)\trace(P_r\circ\dif\nu_t)\\
&=& -(r+1)\sigma_{r+1}-F(\nu_t)\trace(P_r\circ\dif\nu_t),
\end{array}$$
 $$\begin{array}
  {rcl}
  \DIV(P_rX^{\top})
&=& \sum_{ij}((P_r)_{ji}\langle X,
e_i\rangle)_j\\
&=& \sum_{ij}(P_r)_{ji}\delta_{ij}+\sum_{ij}(P_r)_{ji}h_{ij}\langle X, \nu_t\rangle\\
&=& \trace(P_r)-\trace(P_r\circ\dif\nu_t)\langle X, \nu_t\rangle\\
&=& (n-r)\sigma_r-\trace(P_r\circ\dif\nu_t)\langle X, \nu_t\rangle.
\end{array}$$
Thus, the conclusion follows.
\end{proof}
\begin{theorem}\label{thm4.1}
 (First variational formula of $\mathscr{A}_{r, F}$)
\begin{equation}\label{44}
  \mathscr{A}_{r, F}'(t)=-(r+1)\int_M\psi\sigma_{r+1}\dif A_{X_t}.
\end{equation}
\end{theorem}
\begin{proof} We have $(F(\nu_t))'=\langle\grad_{S^n}F, \nu_t'\rangle$,
so by use of Lemma \ref{lemma4.1}, Lemma \ref{newlemma}, (\ref{9}),
(\ref{10}), (\ref{gradf}), (\ref{pr+1}) and Stokes formula, we have
$$\begin{array}
  {rl}
  \mathscr{A}_{r, F}'(t)=&\int_M(F(\nu_t)\sigma_r'+(F(\nu_t))'\sigma_r)\dif
  A_{X_t}+F(\nu_t)\sigma_r\partial_t\dif A_{X_t}\\
  =&\int_M\{F(\nu_t)\DIV(T_{r-1}\grad\psi)+F(\nu_t)\psi\langle T_{r-1}\circ\dif\nu_t,
  \dif\nu_t\rangle+F(\nu_t)\langle\grad\sigma_r, \xi\rangle\\
  &+\langle\sigma_r(\grad_{S^n}F)\circ\nu_t,
  -\grad\psi+\dif\nu_t(\xi)\rangle+F(\nu_t)\sigma_r(-nH\psi+\DIV\xi)\}\dif
  A_{X_t}\\
  =&\int_M\{-\langle\grad (F(\nu_t)),
  T_{r-1}\grad\psi\rangle+F(\nu_t)\psi\langle T_{r-1}\circ\dif\nu_t,
  \dif\nu_t\rangle\\
  &+\langle F(\nu_t)\grad\sigma_r,
  \xi\rangle+\psi\DIV(\sigma_r(\grad_{S^n}F)\circ\nu_t)+\langle\sigma_r\grad (F(\nu_t)), \xi\rangle\\
  &-nH\psi F(\nu_t)\sigma_r
 +F(\nu_t)\sigma_r\DIV\xi\}\dif A_{X_t}\\
 =&\int_M\{-\langle T_{r-1}\grad (F(\nu_t)),
  \grad\psi\rangle+F(\nu_t)\psi\langle T_{r-1}\circ\dif\nu_t,
  \dif\nu_t\rangle\\
  &+\psi\DIV(\sigma_r(\grad_{S^n}F)\circ\nu_t)-nH\psi F(\nu_t)\sigma_r\}\dif A_{X_t}\\
  =&\int_M\psi\{\DIV(\sigma_r(\grad_{S^n}F)\circ\nu_t)+\DIV(T_{r-1}\grad (F(\nu_t)))\\
  &+F(\nu_t)\langle T_{r-1}\circ\dif\nu_t,
  \dif\nu_t\rangle-nHF(\nu_t)\sigma_r\}\dif A_{X_t}\\
  =&\int_M\psi\{\DIV[(\sigma_r+T_{r-1}\circ\dif\nu_t)(\grad_{S^n} F)\circ\nu_t]\\
  &+F(\nu_t)\trace[(T_{r-1}\circ\dif\nu_t+\sigma_rI)\circ \dif\nu_t]\}\dif A_{X_t}\\
  =&\int_M\psi\{\DIV(P_r(\grad_{S^n} F)\circ\nu_t)+F(\nu_t)\trace(P_r\circ \dif\nu_t)\}\dif
  A_{X_t}\\
  =& -(r+1)\int_M\psi\sigma_{r+1}\dif A_{X_t}.
\end{array}$$
\end{proof}

\begin{remark}
 When $F=1$, Lemma 4.1 and Theorem \ref{thm4.1} were proved by R. Reilly
\cite{reilly} (also see \cite{BC}, \cite{CL}).
\end{remark}

  From (\ref{6}), (\ref{11}) and (\ref{44}), we get

\begin{proposition}\label{prop4.4}
(the first variational formula). For all variations of $X$
preserving
 $V$, we have
\begin{equation}\label{48}
\mathscr{A}^\prime_r(t)=\mathscr{F}_{r, F;
\Lambda}'(t)=-\int_M\psi\{(r+1)\sigma_{r+1}-\Lambda\}\dif A_{X_t}.
\end{equation}
\end{proposition}

 Hence we obtain the Euler-Lagrange equation for such a variation
\begin{equation}\label{49}
 (r+1)\sigma_{r+1}-\Lambda=0.
\end{equation}

\begin{theorem}\label{thm4.2}
 (the second variational formula). Let $X: M\to R^{n+1}$ be an
$n$-dimensional closed hypersurface, which satisfies (\ref{49}),
then for all variations of $X$ preserving
 $V$, the second variational formula of  $\mathscr{A}_{r, F}$ at $t=0$ is given by
\begin{equation}\label{50}
\mathscr{A}^{\prime\prime}_r(0)=\mathscr{F}_{r, F; \Lambda}''(0)=
 -(r+1)\int_M\psi\{L_r\psi+\psi\langle T_r\circ\dif\nu, \dif\nu\rangle\}\dif A_X,
\end{equation}
where $\psi$ satisfies
\begin{equation}\label{51}
\int_M \psi\dif A_X=0. \end{equation}
\end{theorem}

\begin{proof}
 Differentiating  (\ref{48}), we get (\ref{50}) by use of (\ref{49}).
\end{proof}

We call $X:M\to R^{n+1}$ to be a stable critical point of
$\mathscr{A}_{r, F}$ for all variations of $X$ preserving  $V$, if
it satisfies (\ref{49}) and $\mathscr{A}^{\prime\prime}_r(0)\geq 0$
for all $\psi$ with condition (\ref{51}).

 \section{Proof of Theorem \ref{thm1.2}}
 Firstly, we prove that if $X(M)$ is, up to translations and homotheties, the Wulff shape, then $X$ is stable.

    From $\dif\phi=(D^2F+F1)\circ\dif x$,
    $\dif\phi$ is perpendicular to $x$. So $\nu=-x$ is the unit inner normal
    vector. We have
   \begin{equation}\label{52}
   \dif\phi=-A_F\circ\dif\nu=\sum_{ijk}A_{jk}h_{ki}\omega_ie_j.
\end{equation}
  On the other hand,
   \begin{equation}\label{53}
   \dif\phi=\sum_i\omega_ie_i,
    \end{equation}
    so we have
\begin{equation}\label{54}
s_{ij}=\sum_kA_{ik}h_{kj}=\delta_{ij}. \end{equation}
    From  this, we easily get $\sigma_r=C_n^r$ and $\sigma_{r+1}=C^{r+1}_n$, thus the Wulff shape satisfies
(\ref{49}) with $\Lambda=(r+1)C^{r+1}_n$. Through a direct
calculation, we
    easily know for Wulff shape,
\begin{equation}\label{55}
\mathscr{A}^{\prime\prime}_r(0)=-(r+1)C^r_{n-1}\int_M[\DIV(A_F\grad\psi)+\psi\langle
A_F\circ\dif\nu,
    \dif\nu\rangle]\dif A_X,
\end{equation}
 and $\psi$ satisfies
\begin{equation}\label{56}
\int_M\psi\dif A_X=0. \end{equation}
 From Palmer \cite{palmer} (also see Winklmann \cite{winklmann}), we know $\mathscr{A}^{\prime\prime}_r(0)\geq 0$, that is,
 the Wulff shape is stable.

 Next, we prove that if $X$ is stable, then up to translations and homotheties, $X(M)$ is the Wulff shape.
 We recall the following lemmas:

\begin{lemma}\label{lemma5.1} (\cite{HeL1}, \cite{HeL2})
     For each
$r=0, 1, \cdots, n-1$, the following integral formulas of Minkowski type hold:
\begin{equation}\label{57}
 \int_M(H^F_rF(\nu) +H^F_{r+1}\langle X, \nu\rangle)\dif A_X=0,\quad r=0, 1,
 \cdots, n-1.
\end{equation}
\end{lemma}
\begin{lemma}\label{lemma5.2}(\cite{HeL1}, \cite{HeL2}, \cite{palmer})
  If $\lambda_1=\lambda_2=\cdots=\lambda_n=\mbox{const}\neq0$, then up to translations and
homotheties,
  $X(M)$ is the Wulff shape.
\end{lemma}

      From Lemma \ref{lemma5.1} and (\ref{51}), we can choose $\psi=\alpha F(\nu)+H^F_{r+1}\langle
      X, \nu\rangle$ as the test function, where $\alpha=\int_MF(\nu) H^F_r\dif A_X/\int_MF(\nu)\dif A_X$. For every smooth function
      $f\colon M\to \mathbb{R}$, and each $r$, we define:
     \begin{equation}\label{62}
      I_r[f]=L_rf+f\langle T_r\circ\dif\nu, \dif\nu\rangle,
      \end{equation}
      Then, we have from (\ref{50})
      \begin{equation}\label{63}
        \mathscr{A}^{\prime\prime}_r(0)=
   -(r+1)\int_M\psi I_r[\psi]\dif A_X.
\end{equation}

\begin{lemma}\label{Ir}
For each $0\leq r\leq n-1$, we have
\begin{equation}\label{64}
I_r[F\circ\nu]=-\langle\grad\sigma_{r+1},
      (\grad_{S^n}F)\circ\nu\rangle+\sigma_1\sigma_{r+1}-(r+2)\sigma_{r+2},
\end{equation}
      and
  \begin{equation}\label{65}
   I_r[\langle X, \nu\rangle]=-\langle\grad\sigma_{r+1}, X^{\top}\rangle-(r+1)\sigma_{r+1}.
\end{equation}
\end{lemma}
\begin{proof}
  From (\ref{gradf}) and (\ref{pr+1}),
  $$\begin{array}
    {rcl}
    I_r[F\circ\nu]&=&\DIV\{T_r\grad(F(\nu))\}+F(\nu)\langle
      T_r\circ\dif\nu, \dif\nu\rangle\\
 &=& \DIV(T_r\circ\dif\nu(\grad_{S^n}F)\circ\nu)+F(\nu)\langle
      T_r\circ\dif\nu, \dif\nu\rangle\\
      &=& \DIV(P_{r+1}(\grad_{S^n}F)\circ\nu)+F(\nu)\trace(P_{r+1}\dif\nu)-\langle\grad\sigma_{r+1}, (\grad_{S^n}F)\circ\nu\rangle\\
      &&-\sigma_{r+1}\{\DIV(P_0(\grad_{S^n}F)\circ\nu)+F(\nu)\trace(P_0\dif\nu)\},
  \end{array}$$
    $$\begin{array}
    {rcl}
    I_r[\langle X, \nu\rangle]&=&\DIV(T_r\grad\langle X, \nu\rangle)+\langle X, \nu\rangle\langle
      T_r\circ\dif\nu, \dif\nu\rangle\\
 &=& \DIV(T_r\circ\dif\nu X^{\top})+\langle X, \nu\rangle\langle
      T_r\circ\dif\nu, \dif\nu\rangle\\
      &=& \DIV(P_{r+1}X^{\top})+\langle X, \nu\rangle\trace(P_{r+1}\dif\nu)-\langle\grad\sigma_{r+1}, X^{\top}\rangle\\
      &&-\sigma_{r+1}\{\DIV(P_0X^{\top})+\langle X,
      \nu\rangle\trace(P_0\dif\nu)\}.
  \end{array}$$
So the conclusions follow from Lemma \ref{newlemma}.
\end{proof}

As $H^F_{r+1}$ is a constant, from (\ref{64}) and (\ref{65}), we
have
\begin{equation}\label{66}
\begin{array}{rcl}
 I_r[\psi] &=& \alpha I_r[F\circ\nu]+H^F_{r+1}I_r[\langle X,
 \nu\rangle]\\
 &=&\alpha(\sigma_1\sigma_{r+1}-(r+2)\sigma_{r+2})-(r+1)H^F_{r+1}\sigma_{r+1}\\
 &=& C^{r+1}_n\{\alpha [nH^F_1H^F_{r+1}-(n-r-1)H^F_{r+2}]-(r+1)(H^F_{r+1})^2\}.
\end{array}
\end{equation}
     Therefore we obtain from Lemma \ref{lemma5.1} (recall
     $H^F_{r+1}$ is constant and $\int_M\psi\dif A_X=0$)
 $$\begin{array}{lcl}
&&\dfrac{1}{r+1}\mathscr{A}^{\prime\prime}_r(0)\\
&=&-\int_M\psi I_r[\psi]dA_X\\
&=&-\int_M\psi C^{r+1}_n\{\alpha
[nH^F_1H^F_{r+1}-(n-r-1)H^F_{r+2}]-(r+1)(H^F_{r+1})^2\}\dif
A_X\\
&=&-\alpha C^{r+1}_n\int_M[\alpha F(\nu)+H^F_{r+1}\langle X, \nu\rangle ][nH^F_1H^F_{r+1}-(n-r-1)H^F_{r+2}]\dif A_X\\
&=&-\alpha^2C_n^{r+1}\int_MF(\nu)[nH^F_1H^F_{r+1}-(n-r-1)H^F_{r+2}]\dif A_X\\
&&-\alpha C^{r+1}_nH^F_{r+1}\int_M\langle X,
\nu\rangle[nH^F_1H^F_{r+1}
-(n-r-1)H^F_{r+2}]\dif A_X\\
&=&-\alpha^2C_n^{r+1}\int_MF(\nu)[nH^F_1H^F_{r+1}-(n-r-1)H^F_{r+2}]\dif A_X\\
&&+\alpha C^{r+1}_nH^F_{r+1}\int_MF(\nu)[nH^F_{r+1}
-(n-r-1)H^F_{r+1}]\dif A_X\\
 &=&-\alpha^2(n-r-1)C_n^{r+1}\int_MF(\nu)(H^F_1H^F_{r+1}-H^F_{r+2})\dif
A_X\vspace*{2mm}\\
 &&-\dfrac{\alpha(r+1)C_n^{r+1}(H^F_{r+1})^2}{\int_MF(\nu)\dif A_X}\{\int_MF(\nu) H^F_1\dif
A_X\int_MF(\nu)\dfrac{H^F_r}{H^F_{r+1}}\dif A_X-(\int_MF(\nu)\dif
A_X)^2\}.
\end{array}$$

As $H^F_{r+1}$ is a constant, it must be positive by the compactness
of $M$. Thus,  by Lemma \ref{lemma3.0}, $H^F_1, \cdots, H^F_r$ are
all positive. So, from \cite{hardy} or \cite{yano}, we have:

(\romannumeral1) for each $0\leq r<n-1$,
\begin{equation}\label{68}
  H^F_1H^F_{r+1}-H^F_{r+2}\geq0,
\end{equation}
with the equality holds if and only if $\lambda_1=\cdots=\lambda_n$,
and

(\romannumeral2) for each $1\leq r\leq n-1$,
\begin{equation}\label{69}
\begin{array}{rcl}
&& \int_MF(\nu) H^F_1\dif
A_X\int_MF(\nu)\dfrac{H^F_r}{H^F_{r+1}}\dif A_X-(\int_MF(\nu)\dif A_X)^2\\
&\geq& \int_MF(\nu) H^F_1\dif A_X\int_MF(\nu)/H^F_1\dif
A_X-(\int_MF(\nu)\dif
A_X)^2\\
&\geq & 0,
\end{array}
\end{equation}
with the equality holds if and only if $\lambda_1=\cdots=\lambda_n$.

From (\ref{68}) and (\ref{69}), we easily obtain that, for each
$0\leq r\leq n-1$,
$$\mathscr{A}^{\prime\prime}_r(0)\leq 0,$$
with the equality holds if and only if $\lambda_1=\cdots=\lambda_n$.
Thus, from Lemma \ref{lemma5.2}, up to translations and homotheties,
  $X(M)$ is the Wulff shape. We complete the proof of Theorem \ref{thm1.2}.

\end{document}